\documentstyle[11pt,amstex,amssymb]{amsart}
\newtheorem{thm}{Theorem}[section]
\newtheorem{prop}[thm]{Proposition}
\newtheorem{lemma}[thm]{Lemma}

\newtheorem{remark}[thm]{Remark}

\newcommand{\proof}{\noindent{\it Proof.}\enspace}

\def\bC{{\Bbb C}}
\def\bR{{\Bbb R}}
\def\cF{{\cal F}}
\def\bZ{{\Bbb Z}}
\def\bN{{\Bbb N}}
\def\Ad{{\rm Ad}\,}
\def\cA{{\cal A}}
\def\cB{{\cal B}}
\def\cE{{\cal E}}
\def\cT{{\cal T}}
\def\ET{\cE\cT}
\def\ETf{\ET^{\rm f}}
\def\cS{{\cal S}}
\def\Tr{{\rm Tr}}
\def\cH{{\cal H}}
\def\SM{S_{\rm M}}
\def\<{\langle}
\def\>{\rangle}
\def\1{{\bf 1}}
\def\ffi{\varphi}
\def\eps{\varepsilon}

\begin{document}

\title[Equilibrium States and Entropy Densities]
{Equilibrium states and their entropy densities \\ in gauge-invariant
$C^*$-systems}
\author[N. Akiho]{Nobuyuki Akiho}
\address{Graduate School of Information Sciences,
Tohoku University, Aoba-ku, Sendai 980-8579, Japan}
\author[F. Hiai]{Fumio Hiai$\,^1$}
\address{Graduate School of Information Sciences,
Tohoku University, Aoba-ku, Sendai 980-8579, Japan}
\author[D. Petz]{D\'enes Petz$\,^2$}
\address{Alfr\'ed R\'enyi Institute of Mathematics, Hungarian 
Academy of Sciences, H-1053 Budapest, Re\'altanoda u.\ 13-15, Hungary}
\thanks{$^1\,$Supported in part by
Japan-Hungary Joint Research Project (JSPS) and by the program
``R\&D support scheme for funding selected IT proposals" of the
Ministry of Public Management, Home Affairs, Posts and
Telecommunications.}
\thanks{$^2\,$Supported in part by MTA-JSPS project (Quantum Probability
and Information Theory) and by OTKA T032662.}

\maketitle

\begin{abstract}
A gauge-invariant $C^*$-system is obtained as the fixed point subalgebra of
the infinite tensor product of full matrix algebras under the tensor product
unitary  action of a compact group. In the paper, thermodynamics is studied in
such systems and the chemical potential theory developed by Araki, Haag,
Kastler and Takesaki is used. As a generalization of  quantum spin system,
the equivalence of the KMS condition, the Gibbs condition and the variational
principle is shown for translation-invariant states. The entropy density of
extremal equilibrium states is also investigated in relation to macroscopic
uniformity.
\end{abstract}

\section*{Introduction}

The rigorous treatment of the statistical mechanics of quantum lattice 
(or spin) systems has been one of the major successes of the $C^*$-algebraic 
approach to quantum physics. The main results are due to many people but
a detailed overview is presented in the monograph \cite{BR}. (Chapter 15 of
\cite{OP} is a concise summary, see also \cite{Se}.) The usual quantum 
spin system is described on the infinite tensor product $C^*$-algebra of full
matrix algebras. Given  an interaction $\Phi$, the local Hamiltonian induces
the local dynamics and the local equilibrium state. The global dynamics and
the global equilibrium states are obtained by a limiting procedure. The
equivalence of the KMS condition, the Gibbs condition and the variational
principle for translation-invariant states is the main essence in the theory;
they were established around 1970 (\cite{Ar0,LR,Ro}). The above mentioned
concepts are used to describe equilibrium states. Recently Araki and Moriya
extended the ideas to fermionic lattice systems \cite{AM}. 

An attempt to extend quantum statistical mechanics from the setting of spin
systems to some approximately finite $C^*$-algebras was made by Kishimoto
\cite{Ki3,Ki4}. Motivated by the chemical potential theory due to Araki, Haag,
Kastler and Takesaki \cite{AHKT}, in our previous paper \cite{HP2} we study
the equivalence of the KMS condition, the Gibbs condition and the variational
principle on approximately finite $C^*$-algebras as a natural extension of
the thermodynamics of one-dimensional quantum  lattice systems. It turned out
that equation (2.8) in the proof of \cite[Theorem 2.2]{HP2} does not hold and
the equivalence formulated in that theorem is recovered here under stronger
conditions. (The error in the proof was pointed out to the authors by
E.~St\o rmer and S.~Neshveyev some years ago.)

A gauge-invariant $C^*$-system is obtained as the fixed point subalgebra of
the infinite tensor product of full matrix algebras under the tensor product
unitary  action of a compact group.  This situation is a typical example of
the chemical  potential theory. The primary aim of the present paper is to
recover the main  results in \cite{HP2} in the restrictive setup of such
gauge-invariant $C^*$-systems. The second aim is to discuss entropy densities
and macroscopic uniformity for extremal equilibrium states in such
$C^*$-systems and to  extend the arguments in \cite{HP1}.

\section{Equilibrium states with chemical potentials }

We begin  by fixing basic notations and terminologies. Let $M_d(\bC)$ be the
algebra of $d\times d$ complex matrices. Let $\cF$ denote a one-dimensional
spin (or UHF) $C^*$-algebra $\bigotimes_{k\in\bZ}\cF_k$ with
$\cF_k:=M_d(\bC)$, and $\theta$ the right shift on $\cF$. Let $G$ be a
separable compact group and $\sigma$ a continuous unitary representation of
$G$ on $\bC^d$ so that a product action $\gamma$ of $G$ on $\cF$ is defined
by $\gamma_g:=\bigotimes_\bZ\Ad\sigma_g$, $g\in G$. Let $\cA:=\cF^\gamma$,
the fixed point subalgebra of $\cF$ for the action $\gamma$ of $G$. For a
finite subset $\Lambda\subset\bZ$ let
$\cF_\Lambda:=\bigotimes_{k\in\Lambda}\cF_k$ and
$\cA_\Lambda:=\cA\cap\cF_\Lambda=\cF_\Lambda^\gamma$, the fixed point
subalgebra for $\gamma|_{\cF_\Lambda}$. Then $\cA$ is an AF $C^*$-algebra
generated by $\{\cA_\Lambda\}_{\Lambda\subset\bZ}$ (\cite[Proposition
2.1]{Pr}). The algebra $\cA$ is called the {\it observable algebra} while
$\cF$ is called the {\it field algebra}. Let $\cS(\cA)$ denote the state
space of $\cA$ and $\cS_\theta(\cA)$ the set of all $\theta$-invariant states
of $\cA$.

An {\it interaction} $\Phi$ is a mapping from the finite subsets of $\bZ$
into $\cA$ such that $\Phi(\emptyset)=0$ and $\Phi(X)=\Phi(X)^*\in\cA_X$ for
each finite $X\subset\bZ$. Given an interaction $\Phi$ and a finite subset
$\Lambda\subset\bZ$, define the {\it local Hamiltonian} $H_\Lambda$ by
$$
H_\Lambda:=\sum_{X\subset\Lambda}\Phi(X),
$$
and the {\it surface energy} $W_\Lambda$ by
$$
W_\Lambda:=\sum\{\Phi(X):
X\cap\Lambda\ne\emptyset,\,X\cap\Lambda^c\ne\emptyset\},
$$
whenever the sum converges in norm.

Throughout the paper we assume that an interaction $\Phi$ is {\it
$\theta$-invariant} and has {\it relatively short range}\,; namely,
$\theta(\Phi(X))=\Phi(X+1)$, where $X+1:=\{k+1:k\in X\}$, for every finite
$X\subset\bZ$ and
$$
|||\Phi|||:=\sum_{X\ni0}{\|\Phi(X)\|\over|X|}<\infty,
$$
where $|X|$ means the cardinality of $X$. Let $\cB(\cA)$ denote the set of all
such interactions, which is a real Banach space with the usual linear
operations and the norm $|||\Phi|||$. Moreover, let $\cB_0(\cA)$ denote the
set of all $\Phi\in\cB(\cA)$ such that
$$
\sum_{X\ni0}\|\Phi(X)\|<\infty\quad\mbox{and}\quad
\sup_{n\ge1}\|W_{[1,n]}\|<\infty.
$$
Then $\cB_0(\cA)$ is a real Banach space with the norm
$$
\|\Phi\|_0:=\sum_{X\ni0}\|\Phi(X)\|+\sup_{n\ge1}\|W_{[1,n]}\|
\ \ (\ge|||\Phi|||).
$$
We define the real Banach space $\cB_0(\cF)$ in a similar manner.

When $\Phi\in\cB_0(\cA)$ we have a strongly continuous one-parameter
automorphism group $\alpha^\Phi$ of $\cF$ such that
$$
\lim_{l,m\to\infty}\|\alpha^\Phi_t(a)-
e^{itH_{[-l,m]}}ae^{-itH_{[-l,m]}}\|=0
$$
for all $a\in\cF$ uniformly for $t$ in finite intervals (see \cite[Theorem
8]{Ki1} and also \cite[6.2.6]{BR}). It is straightforward to see that
$\alpha^\Phi_t\theta=\theta\alpha^\Phi_t$ and
$\alpha^\Phi_t\gamma_g=\gamma_g\alpha^\Phi_t$ for all $t\in\bR$ and $g\in G$
so that $\alpha^\Phi_t(\cA)=\cA$, $t\in\bR$. The sextuple
$(\cF,\cA,G,\alpha^\Phi,\gamma,\theta)$ is a so-called {\it field system} in
the chemical potential theory (\cite{AHKT}, \cite[\S5.4.3]{BR}). The most
general notion of equilibrium states is described by the KMS condition in a
general one-parameter $C^*$-dynamical system (see \cite[\S5.3.1]{BR} for
example). In this paper we consider only $(\alpha^\Phi,\beta)$-KMS states with
$\beta=1$; so we refer to those states as just {\it $\alpha^\Phi$-KMS states}.
The next proposition says that the $\alpha^\Phi$-KMS states are automatically
$\theta$-invariant. This was stated in \cite[Proposition 4.2]{HP2} but the
proof there was given in a wrong way.

\begin{prop}\label{P-1.1}
Let $\Phi\in\cB_0(\cA)$, and let $K(\cA,\Phi)$ denote the set of all
$\alpha^\Phi$-KMS states of $\cA$. Then $K(\cA,\Phi)\subset\cS_\theta(\cA)$,
and $\omega\in K(\cA,\Phi)$ is extremal in $K(\cA,\Phi)$ if and only if
$\omega$ is extremal in $\cS_\theta(\cA)$.
\end{prop}

\proof
The proof below is essentially same as in \cite[\S III]{FVV}. Recall that
the generator of $\alpha^\Phi$ is the closure of the derivation $\delta_0$
with domain $D(\delta_0)=\bigcup_\Lambda\cA_\Lambda$ (over the finite
intervals $\Lambda\subset\bZ$) given by
$$
\delta_0(a):=i\sum_{X\cap\Lambda\ne\emptyset}[\Phi(X),a],
\qquad a\in\cA_\Lambda.
$$
For each $n\in\bN$ let $u_n\in\cF_{[-n,n]}$ be a unitary implementing the
cyclic permutation of $\cF_{[-n,n]}=\bigotimes_{-n}^nM_d(\bC)$, i.e.,
$$
\Ad u_n(a_{-n}\otimes a_{-n+1}\otimes\cdots\otimes a_{n-1}\otimes a_n)
=a_n\otimes a_{-n}\otimes a_{-n+1}\otimes\cdots\otimes a_{n-1}
$$
for $a_k\in M_d(\bC)$. Since $[u_n,\bigotimes_{-n}^n\sigma_g]=0$, we get
$\gamma_g(u_n)=u_n$ for all $g\in G$ so that $u_n\in\cA$. Moreover, since
$\Ad u_n(a)=\theta(a)$ whenever $a\in\cA_{[-n,n-1]}$, it is immediate to see
that $\theta(a)=\lim_{n\to\infty}\Ad u_n(a)$ for all $a\in\cA$. Hence, one
can apply \cite[Corollary II.3]{FVV} (or \cite[5.3.33A]{BR}) to obtain
$K(\cA,\Phi)\subset\cS_\theta(\cA)$, and it suffices to show that
$\sup_{n\ge1}\|\delta_0(u_n)\|<\infty$. This indeed follows because
\begin{eqnarray*}
\|\delta_0(u_n)\|
&=&\Bigg\|\sum_{X\cap[-n,n]\ne\emptyset}[\Phi(X),u_n]\Bigg\| \\
&=&\Bigg\|\sum_{X\cap[-n,n]\ne\emptyset}(\Phi(X)-u_n\Phi(X)u_n^*)\Bigg\| \\
&\le&\Bigg\|\sum_{X\subset[-n,n-1]}(\Phi(X)-\theta(\Phi(X)))\Bigg\|
+\Bigg\|\sum_{X\cap[-n,n]\ne\emptyset\atop X\not\subset[-n,n-1]}
(\Phi(X)-u_n\Phi(X)u_n^*)\Bigg\| \\
&\le&\Bigg\|\sum_{X\subset[-n,n-1]}(\Phi(X)-\Phi(X+1))\Bigg\|
+2\Bigg\|\sum_{X\cap[-n,n]\ne\emptyset\atop X\not\subset[-n,n-1]}
\Phi(X)\Bigg\| \\
&\le&\sum_{X\ni-n}\|\Phi(X)\|+\sum_{X\ni n}\|\Phi(X)\|
+2\sum_{X\ni n}\|\Phi(X)\|
+2\Bigg\|\sum_{X\cap[-n,n]\ne\emptyset\atop X\not\subset[-n,n]}
\Phi(X)\Bigg\| \\
&\le&4\sum_{X\ni0}\|\Phi(X)\|+2\|W_{[-n,n]}\| \\
&\le&4\|\Phi\|_0<\infty.
\end{eqnarray*}

For each $\omega\in\cS_\theta(\cA)$ let $(\pi_\omega,\cH_\omega)$ be the GNS
cyclic representation of $\cA$ associated with $\omega$ and $U_\theta$ be a
unitary implementing $\theta$ so that
$\pi_\omega(\theta(a))=U_\theta\pi_\omega(a)U_\theta^*$ for $a\in\cA$. Since
$(\cA,\theta)$ is asymptotically abelian in the norm sense, i.e.,
$\lim_{|n|\to\infty}\|[a,\theta^n(b)]\|=0$ for all $a,b,\in\cA$, it is well
known \cite[4.3.14]{BR} that
\begin{equation}\label{F-1.1}
\pi_\omega(\cA)'\cap\{U_\theta\}'
\subset\pi_\omega(\cA)'\cap\pi_\omega(\cA)''.
\end{equation}
According to \cite[Lemma 4.7]{TW}, the second assertion is a consequence of this
together with the first assertion (see also \cite[4.3.17 and 5.3.30\,(3)]{BR}
for extremal points of $\cS_\theta(\cA)$ and of $K(\cA,\Phi)$).\qed

\begin{remark}{\rm
Since $(\cA,\theta)$ is asymptotically abelian as mentioned in the above proof,
$\cS_\theta(\cA)$ becomes a simplex. It is also well known that $K(\cA,\Phi)$
is a simplex. These were shown in \cite[\S4]{TW}, where the lattice (or simplex)
structure of state spaces was discussed in a rather general setting. (See also
\cite[4.3.11 and 5.3.30\,(2)]{BR}). Moreover, it is seen from \eqref{F-1.1}
\cite[Lemma 4.7$'$]{TW} that $K(\cA,\Phi)$ is a face of $\cS_\theta(\cA)$.
}\end{remark}

It is known \cite[Lemma 4.1]{HP2} that any tracial state $\phi$ of $\cA$ is
$\theta$-invariant and $\phi$ is extremal if and only if it is multiplicative
in the sense that $\phi(ab)=\phi(a)\phi(b)$ for all $a\in\cA_{[i,j]}$ and
$b\in\cA_{[j+1,k]}$, $i\le j<k$. The $\theta$-invariance of any tracial state
of $\cA$ is a particular case of Proposition \ref{P-1.1} where $\Phi$ is
identically zero. We denote by $\ETf(\cA)$ the set of all faithful and
extremal tracial states of $\cA$. On the other hand, we denote by
$\Xi(G,\sigma)$ the set of all continuous one-parameter subgroups
$t\mapsto\xi_t$ of $G$. Two elements $\xi,\xi'$ in $\Xi(G,\sigma)$ are
identified if there exists $g\in G$ such that
$\Ad\sigma_{g^{-1}\xi_tg}=\Ad\sigma_{\xi_t'}$, $t\in\bR$. In fact, this
defines an equivalence relation and we redefine $\Xi(G,\sigma)$ as the set of
equivalence classes. Then, \cite[Proposition 4.3]{HP2} says

\begin{prop}\label{P-1.3}
There is a bijective correspondence $\phi\leftrightarrow\xi$ between
$\ETf(\cA)$ and $\Xi(G,\sigma)$ under the condition that $\phi$ extends to a
$\gamma_\xi$-KMS state of $\cF$.
\end{prop}

Let $\tau_0$ be the normalized trace on $M_d(\bC)$. Let $\phi$ and $\xi$ be
as in the above proposition. Then there exists a unique selfadjoint
$h\in\cF_{\{0\}}=M_d(\bC)$ such that $\tau_0(e^{-h})=1$ and
$\Ad\sigma_{\xi_t}=\Ad e^{ith}$ for all $t\in\bR$. We call this $h$ the {\it
generator} of $\xi$. Note that $\tau_0(e^{-h}\,\cdot)$ is a unique KMS state
of $M_d(\bC)$ with respect to $\Ad e^{ith}$ and thus
$\hat\phi:=\bigotimes_\bZ\tau_0(e^{-h}\,\cdot)$ is a unique KMS state of
$\cF$ with respect to $\gamma_{\xi_t}=\bigotimes_\bZ\Ad e^{ith}$; so
$\phi=\hat\phi|_\cA$.

Let $\Phi\in\cB_0(\cA)$ and $\xi\in\Xi(G,\sigma)$, and let $\omega$ be an
$\alpha^\Phi$-KMS state of $\cA$. We say that $\xi$ is the {\it chemical
potential} of $\omega$ if there exists an extension $\hat\omega$ of $\omega$
to $\cF$ which is a KMS state with respect to $\alpha^\Phi_t\gamma_{\xi_t}$.
Let $h$ be the generator of $\xi$, and define a $\theta$-invariant interaction
$\Phi^h$ in $\cF$ by
\begin{equation}\label{F-1.2}
\Phi^h(X):=
\begin{cases} \Phi(\{j\})+\theta^j(h) & \text{if $X=\{j\}$,\ $j\in\bZ$}, \\
\Phi(X) & \text{otherwise}. \end{cases}
\end{equation}
Since $\Phi^h\in\cB_0(\cF)$, it generates a one-parameter automorphism group
$\alpha^{\Phi^h}$ on $\cF$. Then, we have
$\alpha^{\Phi^h}_t=\alpha^\Phi_t\gamma_{\xi_t}$, $t\in\bR$, and
$\alpha^\Phi|_\cA=\alpha^{\Phi^h}|_\cA$ (\cite[Lemma 4.4]{HP2}). Due to the
uniqueness of an $\alpha^{\Phi^h}$-KMS state of $\cF$ (\cite{Ar1,Ki2}), we
notice that there is a unique $\alpha^\Phi$-KMS state with chemical
potential $\xi$, which is automatically $\theta$-invariant and faithful. On
the other hand, a consequence of the celebrated chemical potential theory in
\cite[\S II]{AHKT} together with Proposition \ref{P-1.1} is the following: If
$\omega$ is a faithful and extremal $\alpha^\Phi$-KMS state of $\cA$, then
$\omega$ enjoys the chemical potential. A complete conclusion in this
direction will be given in Theorem \ref{T-1.6} below, and Proposition
\ref{P-1.3} is its special case.

To introduce the Gibbs condition, one needs the notion of perturbations of
states of $\cA$. Let $\omega,\psi\in\cS(\cA)$. For each finite interval
$\Lambda\subset\bZ$, the {\it relative entropy} of
$\psi_\Lambda:=\psi|_{\cA_\Lambda}$ with respect to
$\omega_\Lambda:=\omega|_{\cA_\Lambda}$ is given by
$$
S(\psi_\Lambda,\omega_\Lambda)
:=\Tr_\Lambda\biggl({d\psi_\Lambda\over d\Tr_\Lambda}\biggl(
\log{d\psi_\Lambda\over d\Tr_\Lambda}-\log{d\omega_\Lambda\over d\Tr_\Lambda}
\biggr)\biggr).
$$
Here, $\Tr_\Lambda$ denotes the canonical trace on $\cA_\Lambda$ such that
$\Tr_\Lambda(e)=1$ for any minimal projection $e$ in $\cA_\Lambda$. Then the
relative entropy $S(\psi,\omega)$ is defined by
$$
S(\psi,\omega):=\sup_{\Lambda\subset\bZ}S(\psi_\Lambda,\omega_\Lambda)
=\lim_{n\to\infty}S(\psi_{[-n,n]},\omega_{[-n,n]}).
$$
(See \cite{OP} for details on the relative entropy for states of a
$C^*$-algebra.) For each $\omega\in\cS(\cA)$ and $Q=Q^*\in\cA$, since
$\psi\mapsto S(\psi,\omega)+\psi(Q)$ is weakly* lower semicontinuous and
strictly convex on $\cS(\cA)$, the {\it perturbed state} $[\omega^Q]$ by $Q$
is defined as a unique minimizer of this functional (\cite{Do,OP}). Recall
\cite{Ar2,Do} that
\begin{equation}\label{F-1.3}
|S(\psi,\omega)-S(\psi,[\omega^Q])|\le2\|Q\|
\end{equation}
for every $\psi,\omega\in\cS(\cA)$ and $Q=Q^*\in\cA$.

Let $\Phi$ be an interaction in $\cA$ and $\phi$ a tracial state of $\cA$.
For each finite $\Lambda\subset\bZ$, the {\it local Gibbs state}
$\phi^G_\Lambda$ of $\cA_\Lambda$ with respect to $\Phi$ and $\phi$ is
defined by
$$
\phi^G_\Lambda(a):={\phi(e^{-H_\Lambda}a)\over\phi(e^{-H_\Lambda})},
\qquad a\in\cA_\Lambda.
$$
Let $\omega\in\cS(\cA)$ and $(\pi_\omega,\cH_\omega,\Omega_\omega)$ be the
cyclic representation of $\cA$ associated with $\omega$. We say that $\omega$
satisfies the {\it strong Gibbs condition} if $\Omega_\omega$ is separating
for $\pi_\omega(\cA)''$ and if, for each finite $\Lambda\subset\bZ$, there
exists a conditional expectation from $\pi_\omega(\cA)''$ onto
$\pi_\omega(\cA_\Lambda)\vee\pi_\omega(\cA_{\Lambda^c})''$ with respect to
$[\omega^{-W_\Lambda}]\,\tilde{}\,$ and
\begin{equation}\label{F-1.4}
[\omega^{-W_\Lambda}](ab)=\phi^G_\Lambda(a)[\omega^{-W_\Lambda}](b),
\qquad a\in\cA_\Lambda,\ b\in\cA_{\Lambda^c}.
\end{equation}
Here, $[\omega^{-W_\Lambda}]\,\tilde{}\,$ is the normal extension of the
perturbed state $[\omega^{-W_\Lambda}]$ to $\pi_\omega(\cA)''$ (see
\cite[p.~826]{HP2}). Furthermore, we say that $\omega$ satisfies the {\it
weak Gibbs condition} with respect to $\Phi$ and $\phi$ if
$[\omega^{-W_\Lambda}]|_{\cA_\Lambda}=\phi^G_\Lambda$ for any finite
$\Lambda\subset\bZ$.

Now, let $\Phi\in\cB(\cA)$, $\phi\in\ETf(\cA)$ and $\omega\in\cS_\theta(\cA)$.
From now on, for simplicity we write $\cA_n:=\cA_{[1,n]}$, $H_n:=H_{[1,n]}$,
$\phi_n:=\phi|_{\cA_n}$, $\omega_n:=\omega|_{\cA_n}$, etc.\ for each
$n\in\bN$. The {\it mean relative entropy} of $\omega$ with respect to $\phi$
is defined by
$$
\SM(\omega,\phi):=\lim_{n\to\infty}{1\over n}S(\omega_n,\phi_n)
=\sup_{n\ge1}{1\over n}S(\omega_n,\phi_n).
$$
(See \cite[Lemma 3.1]{HP2} for justification of the definition.) Define the
{\it mean energy} $A_\Phi$ of $\Phi$ by
$$
A_\Phi:=\sum_{X\ni0}{\Phi(X)\over|X|}\ \ (\in\cA).
$$
Furthermore, it is known \cite[Theorem 3.5]{HP2} that
$\lim_{n\to\infty}{1\over n}\log\phi(e^{-H_n})$ exists and
$$
\lim_{n\to\infty}{1\over n}\log\phi(e^{-H_n})
=\sup\{-\SM(\omega,\phi)-\omega(A_\Phi):\omega\in\cS_\theta(\cA)\}.
$$
The {\it pressure} of $\Phi$ with respect to $\phi$ is thus defined by
$$
p(\Phi,\phi):=\lim_{n\to\infty}{1\over n}\log\phi(e^{-H_n}).
$$
We have the variational expressions of $p(\Phi,\phi)$ and $\SM(\omega,\phi)$
as follows.

\begin{prop}\label{P-1.4}
Let $\phi\in\ETf(\cA)$. If $\Phi\in\cB(\cA)$, then
\begin{equation}\label{F-1.5}
p(\Phi,\phi)
=\sup\{-\SM(\omega,\phi)-\omega(A_\Phi):\omega\in\cS_\theta(\cA)\}.
\end{equation}
If $\omega\in\cS_\theta(\cA)$, then
\begin{equation}\label{F-1.6}
-\SM(\omega,\phi)=\inf\{p(\Phi,\phi)+\omega(A_\Phi):\Phi\in\cB(\cA)\}.
\end{equation}
\end{prop}

\proof
The expression \eqref{F-1.5} was given in \cite[Theorem 3.5]{HP2} as mentioned
above. We further can transform \eqref{F-1.5} into \eqref{F-1.6} by a simple
duality argument. In fact, for each $\omega\in\cS_\theta(\cA)$ define
$f_\omega\in\cB(\cA)^*$, the dual Banach space of $\cB(\cA)$, by
$f_\omega(\Phi):=-\omega(A_\Phi)$, and set
$\Gamma:=\{f_\omega:\omega\in\cS_\theta(\cA)\}$. Then, it is immediate to see
that $\omega\in\cS_\theta(\cA)\mapsto f_\omega\in\Gamma$ is an affine
homeomorphism in the weak* topologies so that $\Gamma$ is a weakly* compact
convex subset of $\cB(\cA)^*$. Define $F:\cB(\cA)^*\to[0,+\infty]$ by
$$
\begin{cases}
F(f_\omega):=\SM(\omega,\phi) & \text{for $\omega\in\cS_\theta(\cA)$}, \\
F(g):=+\infty & \text{if $g\in\cB(\cA)^*\setminus\Gamma$}.
\end{cases}
$$
Then $F$ is a weakly* lower semicontinuous and convex function on
$\cB(\cA)^*$ (see \cite[Proposition 3.2]{HP2}). Since \eqref{F-1.5} means
that
$$
p(\Phi,\phi)=\sup\{g(\Phi)-F(g):g\in\cB(\cA)^*\},
\quad\Phi\in\cB(\cA),
$$
it follows by duality (see \cite[Proposition I.4.1]{ET} for example) that
$$
F(g)=\sup\{g(\Phi)-p(\Phi,\phi):\Phi\in\cB(\cA)\},
\quad g\in\cB(\cA)^*.
$$
Hence, for every $\omega\in\cS_\theta(\cA)$,
\begin{eqnarray*}
\SM(\omega,\phi)&=&\sup\{f_\omega(\Phi)-p(\Phi,\phi):\Phi\in\cB(\cA)\} \\
&=&-\inf\{p(\Phi,\phi)+\omega(A_\Phi):\Phi\in\cB(\cA)\},
\end{eqnarray*}
giving \eqref{F-1.6}.\qed

\bigskip
We say that $\omega$ satisfies the {\it variational principle} with respect
to $\Phi$ and $\phi$ if
\begin{equation}\label{F-1.7}
p(\Phi,\phi)=-\SM(\omega,\phi)-\omega(A_\Phi).
\end{equation}

With the above definitions in mind we have the next theorem, recovering main
results of \cite{HP2} (Corollary 3.11 and Theorem 4.5) in the special setup of
gauge-invariant $C^*$-systems.

\begin{thm}\label{T-1.5}
Let $\Phi\in\cB_0(\cA)$, $\phi\in\ETf(\cA)$ and $\xi\in\Xi(G,\phi)$ with
$\phi\leftrightarrow\xi$ in the sense of Proposition $\ref{P-1.3}$. Then the
following conditions for $\omega\in\cS(\cA)$ are equivalent:
\begin{itemize}
\item[(i)] $\omega$ is an $\alpha^\Phi$-KMS state with chemical potential
$\xi$;
\item[(ii)] $\omega$ satisfies the strong Gibbs condition with respect to
$\Phi$ and $\phi$;
\item[(iii)] $\omega\in\cS_\theta(\cA)$ and $\omega$ satisfies the weak Gibbs
condition with respect to
$\Phi$ and $\phi$;
\item[(iv)] $\omega\in\cS_\theta(\cA)$ and $\omega$ satisfies the variational
principle with respect to
$\Phi$ and $\phi$.
\end{itemize}

Furthermore, there exists a unique $\omega\in\cS(\cA)$ satisfying one
$($hence all\,$)$ of the above conditions.
\end{thm}

\proof
(i) $\Rightarrow$ (ii). Let $\omega$ be an $\alpha^\Phi$-KMS state with
chemical potential $\xi$ and $(\pi_\omega,\cH_\omega,\Omega_\omega)$ be the
associated cyclic representation of $\cA$. It is well known that
$\Omega_\omega$ is separating for $\pi_\omega(\cA)''$ (see \cite[5.3.9]{BR}
for example). According to the proof of \cite[Theorem 2.2, (i) $\Rightarrow$
(ii)]{HP2}, we see that for any finite $\Lambda\subset\bZ$ there exists a
conditional expectation from $\pi_\omega(\cA)''$ onto
$\pi_\omega(\cA_\Lambda)\vee\pi_\omega(\cA_{\Lambda^c})''$ with respect to
$[\omega^{-W_\Lambda}]\,\tilde{}\,$. (Note that this part of the proof of
\cite[Theorem 2.2, (i) $\Rightarrow$ (ii)]{HP2} is valid.) Moreover, the
proof of \cite[Theorem 4.5]{HP2} shows that \eqref{F-1.4} holds for any
finite $\Lambda\subset\bZ$. Hence we obtain (ii).

(ii) $\Rightarrow$ (iii). The proof of \cite[Theorem 2.2, (ii) $\Rightarrow$
(i)]{HP2} guarantees that (ii) implies $\omega\in K(\cA,\Phi)$. Hence
Proposition \ref{P-1.1} gives the $\theta$-invariance of $\omega$.

(iii) $\Rightarrow$ (iv) is contained in \cite[Proposition 3.9]{HP2} proven in
a more general setting.

(iv) $\Rightarrow$ (i). To prove this as well as the last assertion, it
suffices to show that a state $\omega\in\cS(\cA)$ satisfying (iv) is unique.
First, note that the variational principle \eqref{F-1.7} means that
$\Psi\mapsto-\omega(A_\Psi)$ is a tangent functional to the graph of
$p(\cdot,\phi)$ on $\cB_0(\cA)$ at $\Phi$. Let $h\in M_d(\bC)$ be the
generator of $\xi$ and $\Phi^h$ a $\theta$-invariant interaction in $\cF$
defined by \eqref{F-1.2}. Since $\Phi^h\in\cB_0(\cF)$, there is a unique
$\alpha^{\Phi^h}$-KMS state $\hat\omega$ of $\cF$. Equivalently, there is a
unique $\theta$-invariant state $\hat\omega$ of $\cF$ satisfying the
variational principle with respect to $\Phi^h$, i.e.,
$$
P_\cF(\Phi^h)=s_\cF(\hat\omega)-\hat\omega(A_{\Phi^h}).
$$
Recall here that the pressure $P_\cF(\Psi)$ of $\Psi\in\cB_0(\cF)$ and the
mean entropy $s_\cF(\psi)$ of $\psi\in\cS_\theta(\cF)$ are
$$
P_\cF(\Psi):=\lim_{n\to\infty}{1\over n}\log\Tr_{\cF_n}(e^{-H_n(\Psi)}),
\quad s_\cF(\psi):=\lim_{n\to\infty}{1\over n}S(\psi_n),
$$
where $\Tr_{\cF_n}$ is the usual trace on $\cF_n$ and $H_n(\Psi)$ is the
local Hamiltonian of $\Psi$ inside the interval $[1,n]$. The uniqueness
property above means (see \cite[Proposition I.5.3]{ET} for example) that the
pressure function $P_\cF(\cdot)$ on $\cB_0(\cF)$ is differentiable at
$\Phi^h$. We have (see \cite[(4.11)]{HP2})
\begin{equation}\label{F-1.8}
p(\Phi,\phi)=P_\cF(\Phi^h)-\log d,\qquad\Phi\in\cB_0(\cA).
\end{equation}
By this and \eqref{F-1.2} we obtain
$$p(\Phi+\Psi,\phi)=P_\cF(\Phi^h+\Psi)-\log d,\qquad\Psi\in\cB_0(\cA),
$$
which implies that $\Psi\in\cB_0(\cA)\mapsto p(\Psi,\phi)$ is differentiable
at $\Phi$. Hence the required implication follows.\qed

\bigskip
The next theorem is a right formulation of what we wanted to show in
\cite{HP2}, though in the restricted setup of gauge-invariant $C^*$-systems.

\begin{thm}\label{T-1.6}
If $\Phi\in\cB_0(\cA)$ and $\omega\in\cS(\cA)$, then the following
conditions are equivalent:
\begin{itemize}
\item[(i)] $\omega$ is a faithful and extremal $\alpha^\Phi$-KMS state;
\item[(ii)] $\omega$ is $\alpha^\Phi$-KMS with some chemical potential
$\xi\in\Xi(G,\sigma)$;
\item[(iii)] $\omega$ satisfies the strong Gibbs condition with respect to
$\Phi$ and some $\phi\in\ETf(\cA)$;
\item[(iv)] $\omega\in\cS_\theta(\cA)$ and $\omega$ satisfies the weak Gibbs
condition with respect to $\Phi$ and some $\phi\in\ETf(\cA)$;
\item[(v)] $\omega\in\cS_\theta(\cA)$ and $\omega$ satisfies the variational
principle with respect to $\Phi$ and some $\phi\in\ETf(\cA)$.
\end{itemize}
\end{thm}

\proof
In view of Theorem \ref{T-1.5} we only need to prove the equivalence between
(i) and (ii). (i) $\Rightarrow$ (ii) is a consequence of the chemical
potential theory in \cite[\S II]{AHKT} and Proposition \ref{P-1.1} as
mentioned above (after Proposition \ref{P-1.3}). Conversely, suppose (ii) and
let $\hat\omega$ be a (unique) KMS state of $\cF$ with respect to
$\alpha^\Phi\gamma_{\xi_t}=\alpha^{\Phi^h}$ so that $\omega=\hat\omega|_\cA$.
Since $\hat\omega$ is obviously faithful, so is $\omega$. Moreover, the
extremality of $\omega$ in $\cS_\theta(\cA)$ follows from that of $\hat\omega$
in $\cS_\theta(\cF)$. This may be well known but we sketch the proof for
convenience. Let $(\hat\pi,\hat\cH,\hat\Omega,\hat U_\theta)$ be the cyclic
representation of $\cF$ associated with $\hat\omega$, where $\hat U_\theta$ is
a unitary implementing $\theta$ so that $\hat U_\theta\hat\Omega=\hat\Omega$
and $\hat\pi(\theta(a))=\hat U_\theta\hat\pi(a)\hat U_\theta^*$ for $a\in\cF$.
Then the cyclic representation of $\cA$ associated with $\omega$ is given by
$\cH_\omega:=\overline{\hat\pi(\cA)\hat\Omega}$ and
$\pi_\omega(a):=\hat\pi(a)|_{\cH_\omega}$ for $a\in\cA$ with
$\Omega_\omega:=\hat\Omega$. Let $P:\hat\cH\to\cH_\omega$ be the orthogonal
projection. Since $\hat U_\theta P=P\hat U_\theta$,
$U_\theta^\omega:=\hat U_\theta|_{\cH_\omega}$ is a unitary implementing
$\theta|_\cA$. Let $\hat\sigma$ denote the modular automorphism group of
$\hat\pi(\cF)''$ associated with $\hat\Omega$. Since
$\hat\sigma_t(\hat\pi(a))=\hat\pi(\alpha^\Phi_t(a))\in\hat\pi(\cA)$ for all
$a\in\cA$, there exists the conditional expectation
$E:\hat\pi(\cF)''\to\hat\pi(\cA)''$ with respect to the
state $\<\,\cdot\,\hat\Omega,\hat\Omega\>$ (\cite{Ta}). Notice that $E$ is
$\theta$-covariant, i.e.,
$E(\hat U_\theta x\hat U_\theta^*)=\hat U_\theta E(x)\hat U_\theta^*$ for all
$x\in\hat\pi(\cF)''$. Now, assume that $\omega_1\in\cS_\theta(\cA)$ and
$\omega_1\le\lambda\omega$ for some $\lambda>0$; hence there exists
$T_1\in\pi_\omega(\cA)'$ with $0\le T_1\le\lambda$ such that
$\omega_1(a)=\<T_1\pi_\omega(a)\Omega_\omega,\Omega_\omega\>$ for $a\in\cA$,
and $T_1U_\theta^\omega=U_\theta^\omega T_1$. Define $T:=T_1P+(\1-P)$ on
$\hat\cH$. Then it is easy to check that $0\le T\le\lambda$,
$T\in\hat\pi(\cA)'$ and $T\hat U_\theta=\hat U_\theta T$. Define
$$
\hat\omega_1(a):=\<TE(\hat\pi(a))\hat\Omega,\hat\Omega\>,\qquad a\in\cF,
$$
which is a state of $\cF$ with $\hat\omega_1|_\cA=\omega_1$ and
$\hat\omega_1\le\lambda\hat\omega$. For any $a\in\cF$ we get
$$
\hat\omega_1(\theta(a))
=\<TE(\hat U_\theta\hat\pi(a)\hat U_\theta^*)\hat\Omega,\hat\Omega\>
=\<TE(\hat\pi(a))\hat\Omega,\hat\Omega\>=\hat\omega_1(a)
$$
so that the extremality of $\hat\omega$ implies $\hat\omega_1=\hat\omega$ and
so $\omega_1=\omega$. Hence $\omega$ is extremal in $\cS_\theta(\cA)$ (hence
in $K(\cA,\Phi)$), and (ii) $\Rightarrow$ (i) is shown.\qed

\section{More about variational principle}
\setcounter{equation}{0}

In this section we consider the variational principle for
$\omega\in\cS_\theta(\cA)$ in terms of the mean entropy and the pressure which
are defined by use of canonical traces on local algebras (not with respect to
a tracial state in $\ETf(\cA)$). Let $\nu$ be the restriction of
$\bigotimes_\bZ\tau_0$ to $\cA$, which is an element of $\ETf(\cA)$
corresponding to the trivial chemical potential $\xi=1$. For each $n\in\bN$
the $n$-fold tensor product $\bigotimes_1^n\sigma$ of the unitary
representation $\sigma$ is decomposed as
$$
\bigotimes_1^n\sigma
=m_1\sigma_1\oplus m_2\sigma_2\oplus\cdots\oplus m_{K_n}\sigma_{K_n},
$$
where $\sigma_i\in\widehat G$, $1\le i\le K_n$, are contained in
$\bigotimes_1^n\sigma$ with multiplicities $m_i$. For $1\le i\le K_n$ let
$d_i$ be the dimension of $\sigma_i$. Then, we have
$\sum_{i=1}^{K_n}m_id_i=d^n$ and
\begin{eqnarray}
\cA_n&=&\bigoplus_{i=1}^{K_n}(M_{m_i}(\bC)\otimes\1_{d_i})
\cong\bigoplus_{i=1}^{K_n}M_{m_i}(\bC), \label{F-2.1} \\
\cF_n\cap\cA_n'&=&\bigoplus_{i=1}^{K_n}(\1_{m_i}\otimes M_{d_i}(\bC))
\cong\bigoplus_{i=1}^{K_n}M_{d_i}(\bC). \label{F-2.2}
\end{eqnarray}
The canonical traces $\Tr_{\cA_n}$ on $\cA_n$ and $\Tr_{\cA_n'}$ on
$\cF_n\cap\cA_n'$ are written as
\begin{eqnarray*}
&&\Tr_{\cA_n}\Biggl(\sum_ia_i\otimes\1_{d_i}\Biggr)
=\sum_i\Tr_{m_i}(a_i),\qquad a_i\in M_{m_i}(\bC),\ 1\le i\le K_n, \\
&&\Tr_{\cA_n'}\Biggl(\sum_i\1_{m_i}\otimes b_i\Biggr)
=\sum_i\Tr_{d_i}(b_i),\qquad b_i\in M_{d_i}(\bC),\ 1\le i\le K_n,
\end{eqnarray*}
where $\Tr_m$ denotes the usual trace on $M_m(\bC)$.

\begin{lemma}\label{L-2.1}
$(1)$\enspace If $\omega\in\cS_\theta(\cA)$, then
$\lim_{n\to\infty}{1\over n}S(\omega_n)$ exists and
$$
\lim_{n\to\infty}{1\over n}S(\omega_n)=-\SM(\omega,\nu)+\log d,
$$
where $S(\omega_n)$ is the von Neumann entropy of $\omega_n$ with respect to
$\Tr_{\cA_n}$, i.e.,
\begin{equation}\label{F-2.3}
S(\omega_n):=-\Tr_{\cA_n}\biggl({d\omega_n\over d\Tr_{\cA_n}}
\log{d\omega_n\over d\Tr_{\cA_n}}\biggr)
=-\omega_n\biggl(\log{d\omega_n\over d\Tr_{\cA_n}}\biggr).
\end{equation}

$(2)$\enspace If $\Phi\in\cB(\cA)$, then
$\lim_{n\to\infty}{1\over n}\log\Tr_{\cA_n}(e^{-H_n})$ exists and
$$
\lim_{n\to\infty}{1\over n}\log\Tr_{\cA_n}(e^{-H_n})
=p(\Phi,\nu)+\log d.
$$
\end{lemma}

\proof
(1)\enspace Notice that
$$
S(\omega_n)
=-S(\omega_n,\nu_n)-\omega_n\biggl(\log{d\nu_n\over d\Tr_{\cA_n}}\biggr).
$$
Representing $\cA_n=\bigoplus_{i=1}^{K_n}(M_{m_i}(\bC)\otimes\1_{d_i})$ as in
\eqref{F-2.1}, we have
$$
d^n{d\nu_n\over d\Tr_{\cA_n}}=\sum_{i=1}^{K_n}d_i\1_{m_i}\otimes\1_{d_i},
$$
because
$$
d^n\nu_n\Biggl(\sum_ia_i\otimes\1_{d_i}\Biggr)
=\sum_i\Tr_{m_i}(a_i)\Tr_{d_i}(\1_{d_i})
=\sum_id_i\Tr_{m_i}(a_i)
$$
for $a_i\in M_{m_i}(\bC)$, $1\le i\le K_n$. Therefore,
\begin{equation}\label{F-2.4}
\1_{\cA_n}\le d^n{d\nu_n\over d\Tr_{\cA_n}}
\le\biggl(\max_{1\le i\le K_n}d_i\biggr)\1_{\cA_n}.
\end{equation}
This implies that
$$
0\le\omega_n\biggl(\log{d\nu_n\over d\Tr_{\cA_n}}\biggr)+n\log d
\le\log\biggl(\max_{1\le i\le K_n}d_i\biggr).
$$
As is well known (see a brief explanation in \cite[p.~844]{HP2} for example),
the representation ring of any compact group has polynomial growth; so we have
\begin{equation}\label{F-2.5}
\lim_{n\to\infty}{1\over n}\log\biggl(\max_{1\le i\le K_n}d_i\biggr)=0.
\end{equation}
This implies the desired conclusion.

(2)\enspace By \eqref{F-2.4} we get
$$
\Tr_{\cA_n}(e^{-H_n})\le d^n\nu_n(e^{-H_n})
\le\biggl(\max_{1\le i\le K_n}d_i\biggr)\Tr_{\cA_n}(e^{-H_n}),
$$
implying the result.\qed

\bigskip
In view of the above lemma we define the {\it mean entropy} of
$\omega\in\cS_\theta(\cA)$ by
$$
s_\cA(\omega):=\lim_{n\to\infty}{1\over n}S(\omega_n)
\ \ (=-\SM(\omega,\nu)+\log d),
$$
and the {\it pressure} of $\Phi\in\cB(\cA)$ by
$$
P_\cA(\Phi):=\lim_{n\to\infty}{1\over n}\log\Tr_{\cA_n}(e^{-H_n})
\ \ (=p(\Phi,\nu)+\log d).
$$
The variational expression \eqref{F-1.5} in case of $\phi=\nu$ is rewritten as
$$
P_\cA(\Phi)=\sup\{s_\cA(\omega)-\omega(A_\Phi):\omega\in\cS_\theta(\cA)\}.
$$

\begin{prop}\label{P-2.2}
Let $\Phi\in\cB_0(\cA)$ and $\xi\in\Xi(G,\sigma)$ with the generator $h$.
Assume that $\xi$ is central, i.e., $\xi_t$ belongs to the center of $G$ for
any $t$ $($this is the case if $G$ is abelian\,$)$. Then $\Phi^h$ defined by
\eqref{F-1.2} is an interaction in $\cA$, and $\omega\in\cS_\theta(\cA)$ is
$\alpha^\Phi$-KMS with chemical potential $\xi$ if and only if it satisfies
the variational principle
\begin{equation}\label{F-2.6}
P_\cA(\Phi^h)=s_\cA(\omega)-\omega(A_{\Phi^h}).
\end{equation}
In particular, $\omega$ is $\alpha^\Phi$-KMS with trivial chemical potential
if and only if it satisfies
$$
P_\cA(\Phi)=s_\cA(\omega)-\omega(A_\Phi).
$$
\end{prop}

\proof
The assumption of $\xi$ being central implies that
$\Ad\sigma_g(\sigma_{\xi_t})=\sigma_{\xi_t}$ for all $g\in G$ and $t\in\bR$.
Hence, it is immediate to see that
$\bigotimes_\Lambda e^{-h}=\exp\bigl(-\sum_{j\in\Lambda}\theta^j(h)\bigr)$ is
in $\cA_\Lambda$ for any finite $\Lambda\subset\bZ$ and so the interaction
$\Phi^h$ is in $\cA$. Let $\phi$ be an element of $\ETf(\cA)$ corresponding to
$\xi$ as in Proposition \ref{P-1.3}. We may show that \eqref{F-2.6} is
equivalent to the variational principle \eqref{F-1.7} with respect to $\phi$.
Since $A_{\Phi^h}=A_\Phi+h$, it suffices to prove the following two
expressions:
\begin{equation}\label{F-2.7}
p(\Phi,\phi)=P_\cA(\Phi^h)-\log d
\end{equation}
and for every $\omega\in\cS_\theta(\cA)$
\begin{equation}\label{F-2.8}
-\SM(\omega,\phi)=s_\cA(\omega)-\omega(h)-\log d.
\end{equation}
Let $H_n(\Phi^h)$ be the local Hamiltonian of $\Phi^h$ inside the interval
$[1,n]$. Since
$$
\phi_n(e^{-H_n})
=\nu_n\Biggl(\Biggl(\bigotimes_{j=1}^ne^{-h}\Biggr)e^{-H_n}\Biggr)
=\nu_n(e^{-H_n(\Phi^h)})
=\Tr_{\cA_n}\biggl({d\nu_n\over d\Tr_{\cA_n}}e^{-H_n(\Phi^h)}\biggr),
$$
we obtain \eqref{F-2.7} thanks to \eqref{F-2.4} and \eqref{F-2.5}. On the
other hand, since
\begin{eqnarray}
-S(\omega_n,\phi_n)
&=&S(\omega_n)+\omega_n\biggl(\log{d\phi_n\over d\Tr_{\cA_n}}\biggr)
\nonumber\\
&=&S(\omega_n)+\omega_n\Biggl(\log\Biggl({d\nu_n\over d\Tr_{\cA_n}}
\bigotimes_{j=1}^ne^{-h}\Biggr)\Biggr) \label{F-2.9}\\
&=&S(\omega_n)-n\omega(h)
+\omega_n\biggl(\log{d\nu_n\over d\Tr_{\cA_n}}\biggr), \nonumber
\end{eqnarray}
the expression \eqref{F-2.8} follows.\qed

\section{Entropy densities}
\setcounter{equation}{0}

From now on let $\cF$, $G$, $\sigma$, $\gamma$, $\cA$, $\theta$, etc.\ be
as in the previous sections. Let $\Phi\in\cB_0(\cA)$ be given and
$\alpha^\Phi$ be the associated one-parameter automorphism group. Furthermore,
let $\phi\in\ETf(\cA)$ and the corresponding $\xi\in\Xi(G,\sigma)$ with
generator $h$ be given as in Proposition \ref{P-1.3}; hence $\phi$ extends to
the $\gamma_\xi$-KMS state $\hat\phi$ of $\cF$. For each $n\in\bN$ we then
have the local Gibbs state of $\cA_n$ with respect to $\Phi$ and $\phi$ given
by
$$
\phi_n^G(a):={\phi(e^{-H_n}a)\over\phi(e^{-H_n})},
\qquad a\in\cA_n,
$$
and the local Gibbs state of $\cF_n$ with respect to $\Phi^h$ given by
$$
\hat\phi_n^G(a)
:={\Tr_{\cF_n}(e^{-H_n(\Phi^h)}a)\over\Tr_{\cF_n}(e^{-H_n(\Phi^h)})}.
\qquad a\in\cF_n.
$$
The notation $\hat\phi_n^G$ is justified as follows: Since
$\bigotimes_1^ne^{-h}$ and $e^{-H_n}$ commute (see the proof of
\cite[Proposition 4.3]{HP2}), $\hat\phi_n^G$ is written as 
\begin{equation}\label{F-3.1}
\hat\phi_n^G(a)
={\Tr_{\cF_n}\bigl((\bigotimes_1^ne^{-h})e^{-H_n}a\bigr)
\over\Tr_{\cF_n}\bigl((\bigotimes_1^ne^{-h})e^{-H_n}\bigr)}
={\hat\phi(e^{-H_n}a)\over\hat\phi(e^{-H_n})},
\qquad a\in\cF_n.
\end{equation}
With these notations we have

\begin{thm}\label{T-3.1}
Let $\omega$ be an $\alpha^\Phi$-KMS state of $\cA$ with chemical potential
$\xi$ and $\hat\omega$ be the $\alpha^\Phi\gamma_\xi$-KMS state of $\cF$
extending $\omega$. Then
\begin{eqnarray*}
\SM(\omega,\phi)&=&\lim_{n\to\infty}{1\over n}S(\phi_n^G,\phi_n)
=\lim_{n\to\infty}{1\over n}S(\hat\phi_n^G,\hat\phi_n) \\
&=&\SM(\hat\omega,\hat\phi)
=-s_\cF(\hat\omega)+\hat\omega(h)+\log d
\end{eqnarray*}
and
$$
s_\cF(\hat\omega)=\lim_{n\to\infty}{1\over n}S(\hat\phi_n^G)
=\lim_{n\to\infty}{1\over n}S(\phi_n^G),
$$
where $s_\cF(\hat\omega):=\lim_{n\to\infty}{1\over n}S(\hat\omega_n)$, the
mean entropy of $\hat\omega$. In particular, if $\xi$ is central, then
$s_\cA(\omega)=s_\cF(\hat\omega)$.
\end{thm}

\proof
The following proof of
$\SM(\omega,\phi)=\lim_{n\to\infty}\frac{1}{n}S(\phi^G_n,\phi_n)$ is a slight
modification of \cite[Theorem 2.1]{MvE}. The proof of Theorem \ref{T-1.5}
says that $\Psi\in\cB_0(\cA)\mapsto p(\Psi,\phi)$ is differentiable at $\Phi$
with the tangent functional $\Psi\in\cB_0(\cA)\mapsto-\omega(A_\Psi)$. Hence
we have
\begin{equation}\label{F-3.2}
\frac{d}{d\beta}\biggl|_{\beta=1}p(\beta \Phi,\phi)
=-\omega(A_\Phi).
\end{equation}
Furthermore, we obtain
\begin{equation}\label{F-3.3}
\frac{d}{d\beta}\biggl|_{\beta=1}\frac{1}{n}\log\phi(e^{-H_n(\beta\Phi)})
=\frac{1}{n}\,\frac{\phi(e^{-H_n}(-H_n))}{\phi(e^{-H_n})}
=-\frac{1}{n}\phi^G_n(H_n),
\end{equation}
and as in \cite{MvE}
\begin{equation}\label{F-3.4}
\lim_{n\to\infty}\frac{d}{d\beta}\biggl|_{\beta=1}
\frac{1}{n}\log\phi(e^{-H_n(\beta\Phi)})
=\frac{d}{d\beta}\biggl|_{\beta=1}p(\beta\Phi,\phi).
\end{equation}
Combining \eqref{F-3.2}--\eqref{F-3.4} yields
$\lim_{n\to\infty}\frac{1}{n}\phi^G_n(H_n)=\omega(A_\Phi)$.
Therefore, Theorem \ref{T-1.5} implies
\begin{eqnarray*}
\SM(\omega,\phi)
&=&-p(\Phi,\phi)-\omega(A_\Phi) \\
&=&\lim_{n\to\infty}\frac{1}{n}
\bigl(-\log\phi(e^{-H_n})-\phi^G_n(H_n)\bigr) \\
&=&\lim_{n\to\infty}\frac{1}{n}
\phi^G_n\biggl(\log\frac{d\phi^G_n}{d\phi_n}\biggr) \\
&=&\lim_{n\to\infty}\frac{1}{n}S(\phi^G_n,\phi_n).
\end{eqnarray*}
On the other hand, $\hat\omega$ satisfies the variational principle with
respect to $\Phi^h$, i.e.,
$$
P_\cF(\Phi^h)=s_\cF(\hat\omega)-\hat\omega(A_{\Phi^h}).
$$
Since $A_{\Phi^h}=A_\Phi+h$, this and \eqref{F-1.8} imply
\begin{eqnarray}
\SM(\omega,\phi)
&=&-p(\Phi,\phi)-\omega(A_\Phi) \nonumber\\
&=&-s_\cF(\hat\omega)+\hat\omega(A_\Phi+h)+\log d-\omega(A_\Phi) \nonumber\\
&=&-s_\cF(\hat\omega)+\hat\omega(h)+\log d. \label{F-3.5}
\end{eqnarray}

Since $d\hat\phi_n/d\Tr_{\cF_n}=d^{-n}\bigotimes_{j=1}^ne^{-h}$,
we have
\begin{eqnarray*}
S(\hat\omega_n,\hat\phi_n)
&=&-S(\hat\omega_n)-\hat\omega_n\biggl(
\log\frac{d\hat\phi_n}{d\Tr_{\cF_n}}\biggr) \\
&=&-S(\hat\omega_n)+\hat\omega\Biggl(
\sum_{j=1}^n\theta^j(h)\Biggr)+n\log d \\
&=&-S(\hat\omega_n)+n\hat\omega(h)+n\log d
\end{eqnarray*}
so that
$$
\SM(\hat\omega,\hat\phi)
=-s_\cF(\hat\omega)+\hat\omega(h)+\log d.
$$
Furthermore,
\begin{eqnarray*}
S(\hat\phi^G_n,\hat\phi_n)
&=&-S(\hat\phi^G_n)+\sum_{j=1}^n\hat\phi^G_n(\theta^j(h))+n\log d \\
&=&-S(\hat\phi^G_n)+\sum_{j=1}^n\hat\phi^G_{[1-j,\,n-j]}(h)+n\log d.
\end{eqnarray*}
By \cite{MvE} we have
$s_\cF(\hat\omega)=\lim_{n\to\infty}\frac{1}{n}S(\hat\phi^G_n)$.
The uniqueness of $\alpha^\Phi\gamma_\xi$ ($=\alpha^{\Phi^h}$)-KMS state
implies that $\hat\phi^G_{[-\ell,m]}\to\hat\omega$ weakly* as
$\ell,m\to\infty$. For each $\eps>0$ one can choose $n_0\in\bN$ such that
$\big|\hat\phi^G_{[-\ell,\,m]}(h)-\hat\omega(h)\big|\le\eps$ for all
$\ell,m\ge n_0$. If $n>2n_0$ and $n_0<j\le n-n_0$, then $j-1\ge n_0$ and
$n-j\ge n_0$ so that
$\big|\hat\phi^G_{[1-j,\,n-j]}(h)-\hat\omega(h)\big|\le\eps$. Hence we have
$$
\Bigg|\frac{1}{n}\sum_{j=1}^n\hat\phi^G_{[1-j,\,n-j]}(h)-\hat\omega(h)\Bigg|
\le\frac{4\|h\|n_0}{n}+\eps.
$$
This shows that
$$
\lim_{n\to\infty}\frac{1}{n}\sum_{j=1}^n
\hat\phi^G_{[1-j,\,n-j]}(h)=\hat\omega(h).
$$
Therefore,
$$
\lim_{n\to\infty}\frac{1}{n}S(\hat\phi^G_n,\hat\phi_n)
=-s_\cF(\hat\omega)+\hat\omega(h)+\log d,
$$
and the proof of the first part is completed.

The last assertion follows from \eqref{F-2.6} and \eqref{F-3.5}.
It remains to prove
\begin{equation}\label{F-3.6}
\lim_{n\to\infty}{1\over n}S(\phi_n^G)
=\lim_{n\to\infty}{1\over n}S(\hat\phi_n^G).
\end{equation}
To prove this we give a lemma.

\begin{lemma}\label{L-3.2}
Under \eqref{F-2.1} and \eqref{F-2.2} let
$$
D^0=\sum_{i=1}^{K_n}D^0_i\otimes\1_{d_i}\in\cA_n,
\quad
D'=\sum_{i=1}^{K_n}\1_{m_i}\otimes D'_i\in\cF_n\cap\cA_n'
$$
with positive semidefinite matrices $D^0_i\in M_{m_i}(\bC)$ and
$D'_i\in M_{d_i}(\bC)$ such that $\Tr_{\cF_n}(D^0D')=1$. Then
$D:=D^0D'$ is a density matrix with respect to
$\Tr_{\cF_n}$. If $D|_{\cA_n}$ is the density matrix of
$\Tr_{\cF_n}(D\,\cdot)|_{\cA_n}$ with respect to $\Tr_{\cA_n}$, then
$$
|S(D|_{\cA_n})-S(D)|\le\log\biggl(\max_{1\le i\le K_n}d_i\biggr),
$$
where $S(D)$ is the von Neumann entropy of $D$ with respect to $\Tr_{\cF_n}$
and $S(D|_{\cA_n})$ is that of $D|_{\cA_n}$ with respect to $\Tr_{\cA_n}$
$($see \eqref{F-2.3}$)$.
\end{lemma}

\proof
The first assertion is obvious. Let $E_{\cA_n}$ denote the conditional
expectation from $\cF_n$ onto $\cA_n$ with respect to $\Tr_{\cF_n}$. Notice
that
\begin{eqnarray*}
S(E_{\cA_n}(D))-S(D)
&=&\Tr_{\cF_n}(D\log D-E_{\cA_n}(D)\log E_{\cA_n}(D)) \\
&=&S(D,E_{\cA_n}(D)),
\end{eqnarray*}
the relative entropy of the densities $D$ and $E_{\cA_n}(D)$ in $\cF_n$. Set
$H^0_i:=D^0_i/\Tr_{m_i}(D^0_i)$,
$H'_i:=D'_i/\Tr_{d_i}(D'_i)$ and
$D_i:=H^0_i\otimes H'_i$. The joint convexity of relative
entropy implies
$$
S(D,E_{\cA_n}(D))\le\sum_{i=1}^{K_n}
\Tr_{m_i}(D^0_i)\Tr_{d_i}(D'_i)S(D_i,E_{\cA_n}(D_i)).
$$
Since $E_{\cA_n}(D_i)=d_i^{-1}H^0_i\otimes\1_{d_i}$, we get
\begin{eqnarray*}
&&S(D_i,E_{\cA_n}(D_i)) \\
&&\quad=\Tr_{\cF_n}\Bigl(D_i\Bigl(\log H^0_i\otimes\1_{d_i}
+\1_{m_i}\otimes\log H'_i-\log H^0_i\otimes\1_{d_i}
+(\log d_i)\1_{m_i}\otimes\1_{d_i}\Bigr)\Bigr) \\
&&\quad=\Tr_{d_i}(H'_i\log H'_i)+\log d_i \\
&&\quad\le\log d_i.
\end{eqnarray*}
Therefore,
\begin{equation}\label{F-3.7}
0\le S(E_{\cA_n}(D))-S(D)\le\log\biggl(\max_{1\le i\le K_n}d_i\biggr).
\end{equation}

Next, since for $a=\sum_ia_i\otimes\1_{d_i}\in\cA_n$
\begin{eqnarray*}
\Tr_{\cF_n}(aD')
&=&\Tr_{\cF_n}\Biggl(\sum_{i=1}^{K_n}a_i\otimes D'_i\Biggr)
=\sum_{i=1}^{K_n}\Tr_{m_i}(a_i)\Tr_{d_i}(D'_i) \\
&=&\Tr_{\cF_n}\Biggl(\Biggl(\sum_{i=1}^{K_n}{\Tr_{d_i}(D'_i)\over d_i}
\1_{m_i}\otimes\1_{d_i}\Biggr)a\Biggr),
\end{eqnarray*}
we get
$$
E_{\cA_n}(D')=\sum_{i=1}^{K_n}{\Tr_{d_i}(D'_i)\over d_i}
\1_{m_i}\otimes\1_{d_i}
$$
so that
$$
E_{\cA_n}(D)=D^0E_{\cA_n}(D')
=\sum_{i=1}^{K_n}{\Tr_{d_i}(D'_i)\over d_i}
D^0_i\otimes\1_{d_i}.
$$
Hence we have
\begin{eqnarray*}
&&S(E_{\cA_n}(D)) \\
&&\quad=-\Tr_{\cF_n}\Biggl(\sum_{i=1}^{K_n}{\Tr_{d_i}(D'_i)\over d_i}
D^0_i\otimes\1_{d_i}\Bigl(\log D^0_i\otimes\1_{d_i}
+\Bigl(\log\Tr_{d_i}(D'_i)-\log d_i\Bigr)
\1_{m_i}\otimes\1_{d_i}\Biggr) \\
&&\quad=-\sum_{i=1}^{K_n}\Tr_{d_i}(D'_i)
\Tr_{m_i}(D^0_i\log D^0_i)
-\sum_{i=1}^{K_n}\Tr_{m_i}(D^0_i)\Tr_{d_i}(D'_i)
\Bigl(\log\Tr_{d_i}(D'_i)-\log d_i\Bigr).
\end{eqnarray*}
On the other hand, since $D|_{\cA_n}$ is
$\sum_{i=1}^{K_n}\Tr_{d_i}(D'_i)D^0_i$ as an element of
$\bigoplus_{i=1}^{K_n}M_{m_i}(\bC)$, we have
\begin{eqnarray*}
S(D|_{\cA_n})
&=&-\sum_{i=1}^{K_n}\Tr_{m_i}\Bigl(\Tr_{d_i}(D'_i)D^0_i
\Bigl(\log D^0_i+\log\Tr_{d_i}(D'_i)\Bigr)\Bigr) \\
&=&-\sum_{i=1}^{K_n}\Tr_{d_i}(D'_i)
\Tr_{m_i}(D^0_i\log D^0_i)
-\sum_{i=1}^{K_n}\Tr_{m_i}(D^0_i)
\Tr_{d_i}(D'_i)\log\Tr_{d_i}(D'_i).
\end{eqnarray*}
Therefore,
$$
S(E_{\cA_n}(D))-S(D|_{\cA_n})
=\sum_{i=1}^{K_n}\Tr_{m_i}(D^0_i)\Tr_{d_i}(D'_i)\log d_i
$$
so that
\begin{equation}\label{F-3.8}
0\le S(E_{\cA_n}(D))-S(D|_{\cA_n})
\le\log\biggl(\max_{1\le i\le K_n}d_i\biggr).
\end{equation}
Combining \eqref{F-3.7} and \eqref{F-3.8} gives the conclusion.\qed

\bigskip\noindent
{\it Proof of \eqref{F-3.6}.}\enspace
Let $\hat D_n^G$ be the density of the local Gibbs state $\hat\phi_n^G$ with
respect to $\Tr_{\cF_n}$, which is written as
\begin{equation}\label{F-3.9}
\hat D_n^G={(\bigotimes_1^ne^{-h})e^{-H_n}\over
\Tr_{\cF_n}\bigl((\bigotimes_1^ne^{-h})e^{-H_n}\bigr)}.
\end{equation}
This is obviously of the form of $D$ in Lemma 3.2, i.e., the product of an
element of $\cA_n$ and an element of $\cF_n\cap\cA_n'$. Furthermore, since
$\Tr_{\cF_n}(\hat D_n^G\,\cdot)|_{\cA_n}=\hat\phi_n^G|_{\cA_n}=\phi_n^G$
thanks to \eqref{F-3.1}, it follows that the density of $\phi_n^G$ with
respect to $\Tr_{\cA_n}$ is $\hat D_n^G|_{\cA_n}$ (in the notation of Lemma
\ref{L-3.2}). Hence, Lemma \ref{L-3.2} implies
$$
\big|S(\phi_n^G)-S(\hat\phi_n^G)\big|
\le\log\biggl(\max_{1\le i\le K_n}d_i\biggr)
$$
so that we obtain \eqref{F-3.6} thanks to \eqref{F-2.5}.\qed

\section{Macroscopic uniformity}
\setcounter{equation}{0}

Let $\phi\in\ETf(\cA)$ and $0<\eps<1$. For each $n\in\bN$ and for each state
$\psi$ of $\cA_n$ we define the two quantities
$$
\beta_\eps(\psi):=\min\{\Tr_{\cA_n}(q):
q\in\cA_n\ \mbox{is a projection with}\ \psi(q)\ge1-\eps\},
$$
$$
\beta_\eps(\psi,\phi_n):=\min\{\phi_n(q):
q\in\cA_n\ \mbox{is a projection with}\ \psi(q)\ge1-\eps\}.
$$
For each state $\psi'$ of $\cF_n$ the quantities $\beta_\eps(\psi')$ and
$\beta_\eps(\psi',\hat\phi_n)$ are defined in a similar way with $\cF_n$
instead of $\cA_n$. The aim of this section is to prove

\begin{thm}\label{T-4.1}
Let $\Phi$, $\phi$, $\xi$ and $h$ be as in Theorem $\ref{T-1.5}$, and let
$\omega$ be an $\alpha^\Phi$-KMS state of $\cA$ with chemical potential $\xi$.
Then, for every $0<\eps<1$,
\begin{eqnarray}
-\SM(\omega,\phi)
&=&\lim_{n\to\infty}{1\over n}\log\beta_\eps(\omega_n,\phi_n) \label{F-4.1}\\
&=&\lim_{n\to\infty}{1\over n}\log\beta_\eps(\phi_n^G,\phi_n)
=\lim_{n\to\infty}{1\over n}\log\beta_\eps(\hat\phi_n^G,\hat\phi_n).
\label{F-4.2}
\end{eqnarray}
Moreover, if $\xi$ is central, then for every $0<\eps<1$,
\begin{eqnarray}
s_\cA(\omega)
&=&\lim_{n\to\infty}{1\over n}\log\beta_\eps(\omega_n)
=\lim_{n\to\infty}{1\over n}\log\beta_\eps(\phi_n^G) \label{F-4.3}\\
&=&\lim_{n\to\infty}{1\over n}\log\beta_\eps(\hat\omega_n)
=\lim_{n\to\infty}{1\over n}\log\beta_\eps(\hat\phi_n^G). \label{F-4.4}
\end{eqnarray}
\end{thm}

To prove the theorem, we modify the proofs of \cite[Theorems 3.1 and 3.3]{HP1}.
Let $\omega$ be as in the theorem and
$(\pi_\omega,\cH_\omega,\Omega_\omega)$ be the cyclic representation of $\cA$
associated with $\omega$. For each $n\in\bN$ set
$$
D_n:={d\omega_n\over d\phi_n}\quad\mbox{and}\quad
D_n^G:={d\phi_n^G\over d\phi_n}={e^{-H_n}\over\phi(e^{-H_n})}.
$$

\begin{lemma}\label{L-4.2}
For every $n\in\bN$,
$$
\log D_n^G-\log D_n\le2\|W_n\|.
$$
\end{lemma}

\proof
For every state $\psi$ of $\cA_n$ let $\tilde\psi$ be the state of
$\pi_\omega(\cA_n)$ such that $\psi=\tilde\psi\circ\pi_\omega|_{\cA_n}$; in
particular, let $\tilde\phi_n^G$ be that for $\phi_n^G$. Moreover, let
$\tilde\omega$ be the normal extension of $\omega$ to $\pi_\omega(\cA)''$;
so $\tilde\omega_n=\tilde\omega|_{\pi_\omega(\cA_n)}$. Note (see
\cite[p.~826]{HP2}) that the normal extension $[\omega^{-W_n}]\,\tilde{}\,$
of $[\omega^{-W_n}]$ coincides with the perturbed state
$[\tilde\omega^{-\pi_\omega(W_n)}]$. There exists the conditional expectation
$E_n$ from $\pi_\omega(\cA)''$ onto $\pi_\omega(\cA_n)$ with respect to
$[\omega^{-W_n}]\,\tilde{}\,$ because $\pi_\omega(\cA_n)$ is globally invariant
under the modular automorphism associated with this state. (See the proof of
\cite[Theorem 2.2, (i) $\Rightarrow$ (ii)]{HP2}; this part of the proof of
\cite[Theorem 2.2]{HP2} is valid.) Then, we successively estimate
\begin{eqnarray}
S(\psi,\omega_n)&=&S(\tilde\psi,\tilde\omega_n)
\le S(\tilde\psi\circ E_n,\tilde\omega) \nonumber\\
&\le&S(\tilde\psi\circ E_n,[\omega^{-W_n}]\,\tilde{}\,)+2\|W_n\| \nonumber\\
&=&S(\tilde\psi\circ E_n,\tilde\phi_n^G\circ E_n)+2\|W_n\| \nonumber\\
&=&S(\psi,\phi_n^G)+2\|W_n\|. \label{F-4.5}
\end{eqnarray}
Here, the first inequality is the monotonicity of relative entropy
(\cite[5.12\,(iii)]{OP}) under the restriction of the states of
$\pi_\omega(\cA)''$ to its subalgebra $\pi_\omega(\cA_n)$, and the second is
due to \eqref{F-1.3}. The second equality follows because Theorem \ref{T-1.5}
((ii) or (iii)) gives $[\omega^{-W_n}]\,\tilde{}\,=\tilde\phi_n^G\circ E_n$.
The last equality is seen by applying the monotonicity of relative entropy in
two ways (or by \cite[5.15]{OP}). We now obtain
$$
\psi(\log D_n^G-\log D_n)
=S(\psi,\omega_n)-S(\psi,\phi_n^G)\le2\|W_n\|
$$
for all states $\psi$ of $\cA_n$, which implies the conclusion.\qed

\begin{lemma}\label{L-4.3}
For the densities $D_n$ and $D_n^G$,
$$
\lim_{n\to\infty}{1\over n}\pi_\omega(-\log D_n)
=\lim_{n\to\infty}{1\over n}\pi_\omega(-\log D_n^G)
=-\SM(\omega,\phi)\1\ \ \mbox{strongly}.
$$
\end{lemma}

\proof
Since $\omega$ is extremal in $\cS_\theta(\cA)$, the mean ergodic theorem says
that
$$
\lim_{n\to\infty}{1\over n}
\pi_\omega\Biggl(\sum_{j=1}^n\theta^j(A_\Phi)\Biggr)
=\omega(A_\Phi)\1\ \ \mbox{strongly}.
$$
Since it follows as in \cite{HP1} that
$$
\lim_{n\to\infty}{1\over n}
\Bigg\|\sum_{j=1}^n\theta^j(A_\Phi)-H_n\Bigg\|=0,
$$
we have
\begin{equation}\label{F-4.6}
\lim_{n\to\infty}{1\over n}\pi_\omega(H_n)=\omega(A_\Phi)\1
\ \ \mbox{strongly}.
\end{equation}
Therefore, we obtain the strong convergence
\begin{eqnarray}
{1\over n}\pi_\omega(-\log D_n^G)
&=&{1\over n}\pi_\omega(H_n)+{1\over n}\bigl(\log\phi(e^{-H_n})\bigr)\1
\nonumber\\
&\longrightarrow&\bigl(\omega(A_\Phi)+p(\Phi,\phi)\bigr)\1
=-\SM(\omega,\phi)\1 \label{F-4.7}
\end{eqnarray}
due to the variational principle of $\omega$ in Theorem \ref{T-1.5}.

Next, let $a_n:=-{1\over n}\log D_n$ and
$b_n:=-{1\over n}\log D_n^G+{2\over n}\|W_n\|$; so
$\pi_\omega(b_n)\to-\SM(\omega,\phi)\1$ strongly by what is already shown. We
get $a_n\le b_n$ by Lemma \ref{L-4.2}, and moreover
\begin{eqnarray*}
a_n&=&-{1\over n}\log{d\omega_n\over d\Tr_{\cA_n}}
+{1\over n}\log{d\phi_n\over d\Tr_{\cA_n}}
\ge{1\over n}\log{d\phi_n\over d\Tr_{\cA_n}} \\
&=&-{1\over n}\sum_{j=1}^n\theta^j(h)
+{1\over n}\log{d\nu_n\over d\Tr_{\cA_n}}
\ge-\|h\|-\log d
\end{eqnarray*}
(see \eqref{F-2.9} and \eqref{F-2.4}). Hence $\{b_n-a_n\}$ is uniformly
bounded. Since
\begin{eqnarray*}
\|\pi_\omega(b_n-a_n)\Omega_\omega\|^2
&\le&\Bigl(\sup_m\|b_m-a_m\|\Bigr)\omega(b_n-a_n) \\
&\longrightarrow&\Bigl(\sup_m\|b_m-a_m\|\Bigr)
\bigl(-\SM(\omega,\phi)+\SM(\omega,\phi)\bigr)=0,
\end{eqnarray*}
we have $\pi_\omega(b_n-a_n)\to0$ strongly because $\Omega_\omega$ is
separating for $\pi_\omega(\cA)''$. Hence
$\pi_\omega(a_n)\to-\SM(\omega,\phi)\1$ strongly.\qed

\begin{lemma}\label{L-4.4}
Let $n(1)<n(2)<\cdots$ be positive integers, and let $a_k\in\cA_{n(k)}$ be a
positive contraction for each $k\in\bN$.
\begin{itemize}
\item[(i)] If $\inf_k\omega(a_k)>0$, then
$$
\lim_{k\to\infty}{1\over n(k)}\log\phi_{n(k)}^G(a_k)=0.
$$
\item[(ii)] If $\inf_k\phi_{n(k)}^G(a_k)>0$, then $\inf_k\omega(a_k)>0$.
\item[(iii)] If $\lim_{k\to\infty}\omega(a_k)=1$, then
$\lim_{k\to\infty}\phi_{n(k)}^G(a_k)=1$.
\end{itemize}
The above assertions {\rm(i)--(iii)} hold also for $\cF_{n(k)}$, $\hat\omega$
and $\hat\phi_{n(k)}^G$ instead of $\cA_{n(k)}$, $\omega$ and $\phi_{n(k)}^G$,
respectively.
\end{lemma}

\proof
The last assertion is contained in \cite[Lemma 3.2]{HP1}.

Let
$$
F(s_1,s_2):=s_1\log{s_1\over s_2}+(1-s_1)\log{1-s_1\over1-s_2},
\qquad0\le s_1,s_2\le1.
$$
If the conclusion of (i) does not hold, then one may assume by taking a
subsequence that $\phi_{n(k)}^G(a_k)\le e^{-n(k)\eta}$, $k\in\bN$, for some
$\eta>0$. Using the monotonicity of relative entropy (\cite[5.12\,(iii)]{OP})
applied to the map $\alpha:\bC^2\to\cA_{n(k)}$,
$\alpha(t_1,t_2):=t_1a_k+t_2(\1-a_k)$, we have
\begin{eqnarray*}
S(\omega_{n(k)},\phi_{n(k)}^G)
&\ge&S(\omega_{n(k)}\circ\alpha,\phi_{n(k)}^G\circ\alpha)=
F(\omega_{n(k)}(a_k),\phi_{n(k)}^G(a_k)) \\
&\ge&-\log2-\omega(a_k)\log\phi_{n(k)}^G(a_k)
-(1-\omega(a_k))\log(1-\phi_{n(k)}^G(a_k)) \\
&\ge&-\log2+n(k)\eta\omega(a_k)
\end{eqnarray*}
and hence
$$
\liminf_{k\to\infty}{1\over n(k)}S(\omega_{n(k)},\phi_{n(k)}^G)
\ge\eta\inf_k\omega(a_k)>0.
$$
This contradicts the equality
$$
\lim_{n\to\infty}{1\over n}S(\omega_n,\phi_n^G)
=\SM(\omega,\phi)+\omega(A_\Phi)+p(\Phi,\phi)=0,
$$
which is seen from
$S(\omega_n,\phi_n^G)=S(\omega_n,\phi_n)+\omega(H_n)+\log\phi(e^{-H_n})$ and
\eqref{F-4.6}. Hence (i) follows.

Furthermore, thanks to the monotonicity of relative entropy as above and
\eqref{F-4.5}, we have
\begin{eqnarray*}
F(\phi_{n(k)}^G(a_k),\omega(a_k))
&\le&S(\phi_{n(k)}^G,\omega_{n(k)}) \\
&\le&S(\phi_{n(k)}^G,\phi_{n(k)}^G)+2\|W_{n(k)}\|=2\|W_{n(k)}\|.
\end{eqnarray*}
This shows the boundedness of $F(\phi_{n(k)}^G(a_k),\omega(a_k))$, from which
(ii) and (iii) are easily verified.\qed

\bigskip\noindent
{\it Proof of \eqref{F-4.1}.}\enspace
For each $\delta>0$ and $n\in\bN$, let $p_n$ be the spectral projection of
$-{1\over n}\log D_n$ corresponding to the interval
$(-\SM(\omega,\phi)-\delta,-\SM(\omega,\phi)+\delta)$. Then we have
\begin{equation}\label{F-4.8}
\exp\bigl(n(-\SM(\omega,\phi)-\delta)\bigr)D_np_n
\le p_n\le\exp\bigl(n(-\SM(\omega,\phi)+\delta)\bigr)D_np_n,
\end{equation}
and Lemma \ref{L-4.3} implies that $\pi_\omega(p_n)\to\1$ strongly as
$n\to\infty$. Choose a sequence $n(1)<n(2)<\cdots$ such that
\begin{equation}\label{F-4.9}
\lim_{k\to\infty}{1\over n(k)}\log\beta_\eps(\omega_{n(k)},\phi_{n(k)})
=\liminf_{n\to\infty}{1\over n}\log\beta_\eps(\omega_n,\phi_n).
\end{equation}
For each $k$ choose a projection $q_k\in\cA_{n(k)}$ such that
$\omega(q_k)\ge1-\eps$ and
\begin{equation}\label{F-4.10}
\log\phi_{n(k)}(q_k)\le\log\beta_\eps(\omega_{n(k)},\phi_{n(k)})+1.
\end{equation}
We may assume that $\pi_\omega(q_k)$ converges to some $y\in\pi_\omega(\cA)''$
weakly. Since $\pi_\omega(p_{n(k)}q_k)\to y$ weakly, we get
$$
\lim_{k\to\infty}\omega(p_{n(k)}q_k)=\<y\Omega_\omega,\Omega_\omega\>
=\lim_{k\to\infty}\omega(q_k)\ge1-\eps
$$
and by \eqref{F-4.8}
$$
\phi(q_k)\ge\phi(p_{n(k)}q_k)
\ge\exp\bigl(n(k)(-\SM(\omega,\phi)-\delta)\bigr)\omega(p_{n(k)}q_k).
$$
These give
\begin{equation}\label{F-4.11}
\liminf_{k\to\infty}{1\over n(k)}\log\phi(q_k)
\ge-\SM(\omega,\phi)-\delta.
\end{equation}
Combining \eqref{F-4.9}--\eqref{F-4.11} yields
$$
\liminf_{n\to\infty}{1\over n}\log\beta_\eps(\omega_n,\phi_n)
\ge-\SM(\omega,\phi)-\delta.
$$
On the other hand, we obtain
$$
\limsup_{n\to\infty}{1\over n}\log\beta_\eps(\omega_n,\phi_n)
\le-\SM(\omega,\phi)+\delta,
$$
because by \eqref{F-4.8}
\begin{eqnarray*}
{1\over n}\log\beta_\eps(\omega_n,\phi_n)
&\le&{1\over n}\log\phi(p_n)
\le-\SM(\omega,\phi)+\delta+{1\over n}\log\omega(p_n) \\
&\le&-\SM(\omega,\phi)+\delta
\end{eqnarray*}
if $n$ is so large that $\omega(p_n)\ge1-\eps$. Thus, the proof of
\eqref{F-4.1} is completed.\qed

\bigskip\noindent
{\it Proof of \eqref{F-4.2}.}\enspace
This can be proven by use of (i)--(iii) of Lemma \ref{L-4.4} similarly to the
proof of \cite[Theorem 3.3]{HP1}. Since the proof of the second inequality is
a bit more involved than the first, we only prove the second.

Let $(\hat\pi,\hat\cH,\hat\Omega)$ be the cyclic representation of $\cF$
associated with $\hat\omega$. For each $\delta>0$ and $n\in\bN$, let $p_n$ be
the spectral projection of $-{1\over n}\log D_n^G$ to
$(-\SM(\omega,\phi)-\delta,\SM(\omega,\phi)+\delta)$. Since
${1\over n}\hat\pi(H_n)\to\hat\omega(A_\Phi)\1=\omega(A_\Phi)\1$ and hence
${1\over n}\hat\pi(-\log D_n^G)\to-\SM(\omega,\phi)\1$ strongly as
\eqref{F-4.6} and \eqref{F-4.7}, it follows that $\hat\pi(p_n)\to\1$ strongly
as $n\to\infty$. Furthermore, we have
\begin{equation}\label{F-4.12}
\exp\bigl(n(-\SM(\omega,\phi)-\delta)\bigr)
{e^{-H_n}p_n\over\phi(e^{-H_n})}
\le p_n\le\exp\bigl(n(-\SM(\omega,\phi)+\delta)\bigr)
{e^{-H_n}p_n\over\phi(e^{-H_n})}.
\end{equation}
Choose $n(1)<n(2)<\cdots$ such that
\begin{equation}\label{F-4.13}
\lim_{k\to\infty}{1\over n(k)}
\log\beta_\eps(\hat\phi_{n(k)}^G,\hat\phi_{n(k)})
=\liminf_{n\to\infty}{1\over n}\log\beta_\eps(\hat\phi_n^G,\hat\phi_n).
\end{equation}
For each $k$ there is a projection $q_k\in\cF_{n(k)}$ such that
$\hat\phi_{n(k)}^G(q_k)\ge1-\eps$ and
\begin{equation}\label{F-4.14}
\log\hat\phi_{n(k)}(q_k)
\le\log\beta_\eps(\hat\phi_{n(k)}^G,\hat\phi_{n(k)})+1.
\end{equation}
Here, we may assume that $\hat\pi(q_k)$ converges to some
$y\in\hat\pi(\cF)''$ weakly. Then we obtain
$$
\lim_{k\to\infty}\hat\omega(p_{n(k)}q_kp_{n(k)})
=\<y\hat\Omega,\hat\Omega\>=\lim_{k\to\infty}\hat\omega(q_k)>0
$$
by Lemma \ref{L-4.4}\,(ii) (for $\hat\omega$ and $\hat\phi_{n(k)}^G$ with
$a_k=q_k$), and hence
\begin{equation}\label{F-4.15}
\lim_{k\to\infty}{1\over n(k)}
\log\hat\phi_{n(k)}^G(p_{n(k)}q_kp_{n(k)})=0
\end{equation}
by Lemma \ref{L-4.4}\,(i) (for $\hat\omega$ and $\hat\phi_{n(k)}^G$ with
$a_k=p_{n(k)}q_kp_{n(k)}$). Furthermore, since $p_n$ commutes with $e^{-H_n}$
and $\bigotimes_1^ne^{-h}$, we obtain
\begin{eqnarray*}
\hat\phi_{n(k)}(q_k)
&=&d^{-n(k)}\Tr_{\cF_{n(k)}}\Biggl(\Biggl(\bigotimes_1^{n(k)}e^{-h}
\Biggr)q_k\Biggr) \\
&\ge&d^{-n(k)}\Tr_{\cF_{n(k)}}\Biggl(\Biggl(\bigotimes_1^{n(k)}e^{-h}
\Biggr)p_{n(k)}q_k\Biggr) \\
&\ge&\exp\bigl(n(k)(-\SM(\omega,\phi)-\delta)\bigr)
{d^{-n(k)}\Tr_{\cF_{n(k)}}\bigl(\bigl(\bigotimes_1^{n(k)}e^{-h}
\bigr)e^{-H_{n(k)}}p_{n(k)}q_k\bigr)\over\phi(e^{-H_{n(k)}})} \\
&=&\exp\bigl(n(k)(-\SM(\omega,\phi)-\delta)\bigr)
{\hat\phi(e^{-H_{n(k)}}p_{n(k)}q_kp_{n(k)})\over\phi(e^{-H_{n(k)}})} \\
&=&\exp\bigl(n(k)(-\SM(\omega,\phi)-\delta)\bigr)
\hat\phi_{n(k)}^G(p_{n(k)}q_kp_{n(k)})
\end{eqnarray*}
using \eqref{F-4.12} and \eqref{F-3.1}. This together with
\eqref{F-4.13}--\eqref{F-4.15} yields
$$
\liminf_{n\to\infty}{1\over n}\log\beta_\eps(\hat\phi_n^G,\hat\phi_n)
\ge-\SM(\omega,\phi)-\delta.
$$
On the other hand, since $\hat\phi_n^G(p_n)\to1$ by Lemma \ref{L-4.4}\,(iii)
(for $\hat\omega$ and $\hat\phi_n^G$), we have $\hat\phi_n^G(p_n)\ge1-\eps$
for large $n$, and for such $n$
$$
{1\over n}\log\beta_\eps(\hat\phi_n^G,\hat\phi_n)
\le{1\over n}\log\hat\phi_n(p_n)\le-\SM(\omega,\phi)+\delta
$$
thanks to \eqref{F-4.12}. Therefore,
$$
\limsup_{n\to\infty}{1\over n}\log\beta_\eps(\hat\phi_n^G,\hat\phi)
\le-\SM(\omega,\phi)+\delta,
$$
completing the proof of \eqref{F-4.2}.\qed

\bigskip\noindent
{\it Proof of \eqref{F-4.3} and \eqref{F-4.4}.}\enspace
Assume that $\xi$ is central. Since $s_\cA(\omega)=s_\cF(\hat\omega)$ by
Theorem \ref{T-3.1}, the assertion \eqref{F-4.4} is contained in
\cite[Theorem 3.3]{HP1}. To prove \eqref{F-4.3}, we first assume that $\xi$ is
trivial. Then, by Lemma \ref{L-2.1}\,(1) and \eqref{F-4.1} (in case of
$\phi=\nu$) we have
\begin{eqnarray*}
s_\cA(\omega)&=&-\SM(\omega,\nu)+\log d \\
&=&\lim_{n\to\infty}{1\over n}\log\beta_\eps(\omega_n,\nu_n)+\log d \\
&=&\lim_{n\to\infty}{1\over n}\log\beta_\eps(\omega_n).
\end{eqnarray*}
The latter equality in the above is readily verified from \eqref{F-2.4} and
\eqref{F-2.5}. The other equality in \eqref{F-4.3} when $\phi=\nu$ is
similarly shown from the first equality in \eqref{F-4.2}. When $\xi$ is not
trivial, we consider $\Phi^h$ belonging to $\cB_0(\cA)$ instead of $\Phi$.
Note that $\omega$ is an $\alpha^{\Phi^h}$-KMS state with trivial chemical
potential and $\phi_n^G$ is the local Gibbs state with respect to
$\Phi^h$ and $\nu$. Hence, the above special case gives the conclusion.\qed

\section{Remarks and problems}
\setcounter{equation}{0}

Some problems as well as related known results are in order.

\medskip\noindent
{\bf 5.1.}\enspace
It is known \cite{Ha,Pr} that the weak*-closure of $\ETf(\cA)$ coincides with
the set $\ET(\cA)$ of all extremal tracial states of $\cA$ as far as $G$ is
a compact connected Lie group. For $\Phi\in\cB_0(\cA)$ let $\cE K(\cA,\Phi)$
denote the set of all extremal $\alpha^\Phi$-KMS states of $\cA$ (see
Proposition \ref{P-1.1}) and $\cE K^{\rm f}(\cA,\Phi)$ the set of all
faithful $\omega\in\cE K(\cA,\Phi)$. Theorems \ref{T-1.5} and \ref{T-1.6} say
that there is a bijective correspondence $\phi\leftrightarrow\omega$ between
$\ETf(\cA)$ and $\cE K^{\rm f}(\cA,\Phi)$. We further know (see
\cite[Theorem 4.6]{HP2}) that the correspondence $\phi\mapsto\omega$ is a
weak*-homeomorphism from $\ETf(\cA)$ onto $\cE K^{\rm f}(\cA,\Phi)$. Upon
these considerations we are interested in the following problems:
\begin{itemize}
\item[(1)] Does the weak*-closure of
$\cE K^{\rm f}(\cA,\Phi)$ coincide with $\cE K(\cA,\Phi)$ (as far as $G$ is
a compact connected Lie group)?
\item[(2)] Does the above $\phi\mapsto\omega$ extend to a weak*-homeomorphism
from $\ET(\cA)$ onto $\cE K(\cA,\Phi)$?
\end{itemize}

\medskip\noindent
{\bf 5.2.}\enspace
In the situation of Theorem \ref{T-3.1} it seems that the equality
$s_\cA(\omega)=s_\cF(\hat\omega)$ holds without the assumption of $\xi$ being
central. This is equivalent to the equality
$s_\cA(\omega)=\lim_{n\to\infty}{1\over n}S(\phi_n^G)$, which is the only
missing point in Theorem \ref{T-3.1}.

\medskip\noindent
{\bf 5.3.}\enspace
The equality $-\SM(\omega,\phi)=\lim_{n\to\infty}
{1\over n}\beta_\eps(\hat\omega_n,\hat\phi_n)$ is missing in Theorem 4.1,
which is equivalent to
\begin{equation}\label{F-5.1}
-\SM(\hat\omega,\hat\phi)=\lim_{n\to\infty}
{1\over n}\log\beta_\eps(\hat\omega_n,\hat\phi_n)
\end{equation}
due to Theorem \ref{T-3.1}. Note that $\hat\phi$ is a product state of $\cF$
and $\hat\omega$ is completely ergodic, i.e., extremal for all $\theta^n$,
$n\ge1$. Thus, the equality \eqref{F-5.1} is an old open problem from the
viewpoint of quantum hypothesis testing in \cite{HP0}, where the weaker result
was proven:
\begin{eqnarray*}
-\SM(\hat\omega,\hat\phi)&\ge&\limsup_{n\to\infty}
{1\over n}\log\beta_\eps(\hat\omega_n,\hat\phi_n), \\
-{1\over1-\eps}\SM(\hat\omega,\hat\phi)&\le&\liminf_{n\to\infty}
{1\over n}\log\beta_\eps(\hat\omega_n,\hat\phi_n).
\end{eqnarray*}
In this connection, it is worthwhile to note that T.~Ogawa and H.~Nagaoka
established in \cite{ON} the equality
$$
-S(\ffi,\psi)=\lim_{n\to\infty}{1\over n}
\log\beta_\eps(\ffi_n,\psi_n)
$$
when $\ffi,\psi$ are states of $M_d(\bC)$ and $\ffi_n,\psi_n$ are the
$n$-fold tensor products of $\ffi,\psi$. The problem of macroscopic uniformity
for states of spin $C^*$-algebras was completely solved in a recent paper by
I.~Bjelakovi\'c et al.\ as follows: If $\ffi$ is an extremal
translation-invariant state of the $\nu$-dimensional spin algebra
$\bigotimes_{\bZ^\nu}M_d(\bC)$, then
$$
s(\ffi)=\lim_{\Lambda\to\bZ^\nu}{1\over|\Lambda|}\log\beta_\eps(\ffi)
$$
for any $0<\eps<1$. See \cite{BKSS} for details.

\medskip\noindent
{\bf 5.4.}\enspace
Although many arguments in this paper as well as in \cite{HP2} work also in
gauge-invariant $C^*$-systems over the multi-dimensional lattice $\bZ^\nu$,
some difficulties arise when we would extend our whole arguments to the
multi-dimensional case. For instance, it does not seem that Proposition
\ref{P-1.1} holds in multi-dimensional gauge-invariant $C^*$-systems. The
proposition is crucial when we use the chemical potential theory as in the
proof of Theorem \ref{T-1.6}. Moreover, the assumption of uniformly bounded
surface energies is sometimes useful in our discussions. In the
multi-dimensional case, the assumption is obviously too strong and, if it is
not assumed, the non-uniqueness of KMS states (or the phase transition) can
occur. Indeed, the uniqueness of $\alpha^{\Phi^h}$-KMS state of $\cF$ is
essential in the proof of Theorem \ref{T-1.5}. Consequently, some new ideas
must be needed to extend the theory to the multi-dimensional setting.

\bigskip\noindent
{\bf Acknowledgments.}
The authors are grateful to Professors E.~St\o rmer and S.~Neshveyev who
pointed out a mistake in our previous paper \cite{HP2} in 2000, and also thank
the referees for their useful suggestions.


\begin{thebibliography}{99}

\bibitem{Ar0}
H. Araki,
On the equivalence of the KMS condition and the variational principle for
quantum lattice systems,
{\it Comm. Math. Phys.} {\bf 38} (1974), 1--10.

\bibitem{Ar1}
H. Araki,
On uniqueness of KMS states of one-dimensional quantum lattice systems,
{\it Comm. Math. Phys.} {\bf 44} (1975), 1--7.

\bibitem{Ar2}
H. Araki,
Relative entropy for states of von Neumann algebras II,
{\it Publ. Res. Inst. Math. Sci.} {\bf 13} (1977), 173--192.

\bibitem{AHKT}
H. Araki, R. Haag, D. Kastler and M. Takesaki,
Extension of KMS states and chemical potential,
{\it Comm. Math. Phys.} {\bf 53} (1977), 97--134.

\bibitem{AM}
H. Araki and H. Moriya,
Equilibrium statistical mechanics of Fermion lattice systems,
{\it Rev. Math. Phys.} {\bf 15} (2003), 93--198.

\bibitem{BKSS} 
I. Bjelakovi\'c, T. Kr\"uger, R. Siegmund-Schultze and A. Szko\l a,
The Shannon-McMillan theorem for ergodic quantum lattice systems,
{\it Invent. Math}. {\bf 155} (2004), 203--222.

\bibitem{BR}
O. Bratteli and D. W. Robinson,
{\it Operator Algebras and Quantum Statistical Mechanics 1, 2},
2nd edition, Springer-Verlag, 2002.

\bibitem{Do}
M. J. Donald,
Relative hamiltonians which are not bounded from above,
{\it J. Funct. Anal.} {\bf 91} (1990), 143--173.

\bibitem{ET}
I. Ekeland and R.Temam,
{\it Convex analysis and variational problems},
Studies in Mathematics and its Applications, Vol. 1,
North-Holland, Amsterdam-Oxford, 1976.

\bibitem{FVV}
M. Fannes, P. Vanheuverzwijn and A. Verbeure,
Quantum energy-entropy inequalities: a new method for proving the absence
of symmetry breaking,
{\it J. Math. Phys.} {\bf 25} (1984), 76--78.

\bibitem{Ha}
D. Handelman,
Extending traces on fixed point $C^\ast$ algebras under Xerox product type
actions of compact Lie groups,
{\it J. Funct. Anal.} {\bf 72} (1987), 44--57.

\bibitem{HP0}
F. Hiai and D. Petz,
The proper formula for relative entropy and its asymptotics in quantum
probability,
{\it Comm. Math. Phys.} {\bf 143} (1991), 99--114.

\bibitem{HP1}
F. Hiai and D. Petz,
Entropy densities for Gibbs states of quantum spin systems,
{\it Rev. Math. Phys.} {\bf 5} (1993), 693--712.

\bibitem{HP2}
F. Hiai and D. Petz,
Quantum mechanics in AF $C^*$-systems,
{\it Rev. Math. Phys.} {\bf 8} (1996), 819--859.

\bibitem{Ki1}
A. Kishimoto,
Dissipations and derivations,
{\it Comm. Math. Phys.} {\bf 47} (1976), 25--32.

\bibitem{Ki2}
A. Kishimoto,
On uniqueness of KMS states of one-dimensional quantum lattice systems,
{\it Comm. Math. Phys.} {\bf 47} (1976), 167--170.

\bibitem{Ki3}
A. Kishimoto,
Equilibrium states of a semi-quantum lattice system,
{\it Rep. Math. Phys.} {\bf 12} (1977), 341--374.

\bibitem{Ki4}
A. Kishimoto,
Variational principle for quasi-local algebras over the lattice,
{\it Ann. Inst. H. Poincar\'e Phys. Th\'eor.} {\bf 30} (1979), 51--59.

\bibitem{LR}
O. E. Lanford III and D. W. Robinson,
Statistical mechanics of quantum spin systems.\ III,
{\it Comm. Math. Phys.} {\bf 9} (1968), 327--338.

\bibitem{MvE}
H. Moriya and A. van Enter,
On thermodynamic limits of entropy densities,
{\it Lett. Math. Phys.} {\bf 45} (1998), 323--330.

\bibitem{ON}
T. Ogawa and H. Nagaoka,
Strong converse and Stein's lemma in quantum hypothesis testing,
{\it IEEE Trans. Inform. Theory} {\bf 46} (2000), 2428--2433.

\bibitem{OP}
M. Ohya and D. Petz,
{\it Quantum Entropy and Its Use}, Springer-Verlag, 1993; 2nd edition, 2004. 

\bibitem{Pr}
G. Price,
Extremal traces on some group-invariant $C^*$-algebras,
{\it J. Funct. Anal.} {\bf 49} (1982), 145--151.

\bibitem{Ro}
D. W. Robinson,
Statistical mechanics of quantum spin system.\ II,
{\it Comm. Math. Phys.} {\bf 7} (1968), 337--348.

\bibitem{Se}
G. L. Sewell,
{\it Quantum Theory of Collective Phenomena},
Clarendon Press, New York, 1986.

\bibitem{Ta}
M. Takesaki,
Conditional expectations in von Neumann algebras,
{\it J. Funct. Anal.} {\bf 9} (1972), 306--321.

\bibitem{TW}
M. Takesaki and M. Winnink,
Local normality in quantum statistical mechanics,
{\it Comm. Math. Phys.} {\bf 30} (1973), 129--152.

\end{thebibliography}
\end{document}